\input amstex
\documentstyle{amsppt}
\NoBlackBoxes

\pagewidth{32pc} \pageheight{46pc} \topskip=\normalbaselineskip

\define\CP{\Bbb C \Bbb P}

\define \C{\Bbb C}

\define \Q{\Bbb Q}

\topmatter
\title Diffeomorphisms, Isotopies and Braid Monodromy Factorizations of 
Plane Cuspidal Curves\endtitle

\rightheadtext{ Diffeomorphisms, Isotopies 
and Braid Monodromy Factorizations} 

\author V. Kharlamov and Vik. S. Kulikov \endauthor
\abstract In this paper we prove that there is an infinite sequence of pairs 
of plane cuspidal curves $C_{m,1}$
and $C_{m,2}$,  
such that the pairs $(\CP^2,C_{m,1})$ and 
$(\CP^2,C_{m,2})$ are 
diffeomorphic, but $C_{m,1}$
and $C_{m,2}$ have
non-equivalent braid monodromy factorizations. These curves give rise to 
the negative solutions of "Dif=Def" and "Dif=Iso" problems 
for plane irreducible cuspidal curves.
In our examples,
$C_{m,1}$ and $C_{m,2}$ are complex conjugated.\endabstract
\thanks \flushpar Partially supported by INTAS-OPEN-97-2072, 
NWO-RFBR-047-008-005, and
RFBR  99-01-01133. This work was done during the stay of the second
author in Strasbourg university.\endthanks
\endtopmatter

\document

\baselineskip 20pt 
Up to our knowledge the following natural questions were open till now: 
does an existence of a diffeomorphism $\CP^2\to \CP^2$ 
transforming an algebraic curve $C_1$ in an algebraic curve $C_2$ imply 
the existence of an equisingular deformation ("Dif=Def") or of a 
diffeotopie ("Dif=Iso") between them? We show that the response to the 
both questions is in the negative. In our examples the curves are cuspidal and irreducible (note that for nodal curves the responses are in the positive, 
see \cite{H}).

Let $C\subset \CP^2$ be a plane algebraic curve of degree $\deg C=d$.  
The  
embedding  $C\subset \CP^2$ is determined, up to diffeotopy, 
by the so-called {\it braid monodromy} 
of $C$ (see \cite{KT}). The latter depends on a choice 
of generic homogeneous coordinates in $\CP^2$ and can be seen as a 
factorization of the {\it full twist} $\Delta _{d}^2$ in the normal
semi-group $B^{+}_{d}$ of quasi-positive braids on $d$ 
strings ($\Delta _{d}^2$ is the so-called Garside element;  
$\Delta _{d}^2=(X_1\cdot ...\cdot X_{d-1})^{d}$ in standard 
generators $X_1,\dots , X_{d-1}$ of $B^{+}_{d}$). If the 
only singularities of $C$ are ordinary cusps and nodes ({\it  
cuspidal curve}), then this factorization can be written as follows
$$
\Delta _{d}^2=\prod_{i}\rho_i^{-1}\sigma_1^{s_i}\rho_i, \,\,
s_i\in (1,2,3), \rho_i\in B_d,
\tag 1
$$
where $s_i=1,2,3$ correspond to branch points, nodes 
and cusps respectively, and $B_{d}$ is the braid group.

Let
$$
h=g_1\cdot ...\cdot g_r 
\tag 2
$$
be such a factorization in $B^{+}_{d}$. The transformations which 
replace in (2) two neighboring factors 
$g_i\cdot g_{i+1}$ by $(g_ig_{i+1}g_i^{-1})\cdot g_i $
or
$g_{i+1}(g_{i+1}^{-1}g_ig_{i+1})$
are called {\it Hurwitz moves}.

For $z\in B_{d}$, we denote by
$$h_z=z^{-1}g_1z\cdot z^{-1}g_2z\cdot ...\cdot z^{-1}g_rz$$
and say that $h_z$ is obtained from 
(2) by simultaneous conjugation by $z$. Two factorizations are 
called {\it Hurwitz and conjugation equivalent} if one can be obtained 
from the other by a finite sequence of Hurwitz moves followed by a 
simultaneous conjugation. If two factorizations are Hurwitz and conjugation 
equivalent then they are said to be of the  same 
{\it braid factorization type}.   As is known, any two factorizations 
of the form (1) corresponding to the same  algebraic curve $C\subset \CP^2$ 
are of 
the same braid factorization type (see, for example, \cite{MT}). 

The notion of braid monodromy factorization extends word by word from 
the algebraic case to the case of so called 
{\it semi-algebraic} or {\it Hurwitz curves} 
(details can be found, for example, in 
\cite{KT}).

If $C_1, C_2 \subset \CP^2$ are two projective cuspidal plane curves
of the same braid factorization type, then (\cite{KT}, Theorem 7.1) 
there exists 
a diffeotopy
$F_t : \CP^2 \to \CP^2$ such that $F_1 (C_1)=C_2$.
\footnote{ Actually,  
in \cite{KT} it is proved the existence of symplectic isotopy. 
It is sufficient to rescale the isotopy constructed there.}

The aim of this note is to prove the following theorems.

\proclaim{Theorem 1} There  
are
two sequences of plane irreducible 
cuspidal curves, $\{C_{m,1}\}$ and $\{C_{m,2}\}$, $m\ge 5$, such that
the pairs $(\CP^2,C_{m,1})$ and $(\CP^2,C_{m,2})$
are diffeomorphic, but 
$C_{m,1}$ and $C_{m,2}$ are not  
isotopic and have different braid  
factorization types.
\endproclaim

\proclaim{Theorem 2} 
Let $X$ be a   
surface of general type with ample
canonical class $K$. 
Suppose that there is 
no  
homeomorphism $h$ of $X$ such that $h^*[K]=-[K]$,
$[K]\in H^2(X;\Q)$.
Then the moduli space of $X$ consists of at least 2 
connected components corresponding to 
$X$ and $\overline{X}$ (the bar states for reversing of complex structure, 
$J \mapsto -J$), and for any $m\geq 5$ 
these two connected components are distinguished by 
the braid  
factorization
types of the branch curves of generic coverings 
$f_m:X\to \CP^2$ and $\overline{f}_m:\overline{X}\to \CP^2$ 
given 
by $mK$ and $m\overline K$, respectively.
\endproclaim

These two theorems will be proved simultaneously. Note also that the curves 
from Theorem 1 constructed below are complex conjugated and belong 
to different components of the space of curves of given degree with given 
types of singularities. 

\demo{Proof} 
By Bombieri theorem, the map $X\to\CP^{r_m}$,
$r_m=\dim H^0(X,mK)-1$, given by $mK$ is an imbedding
if $m\ge 5$. Let $m\ge 5$ and denote
by $X_m$ the image of $X$ under  
the imbedding in $\CP^{r_m}$ given by $mK$, by
$pr_m: \CP ^{r_m}\to \CP ^2$ a linear projection generic with respect to 
$X_m$, by $f_m=pr_m|_{X_m}$ the restriction of $pr_m$ to $X_m$, and 
by $C_{m,1}\subset \CP^2$ the branch curve of $f_m$. 

We identify  $X_m$ and $\overline{X}_m$ as sets. Thus, 
the branched covering $f_m:X_m\to \CP^2$ can be considered as 
well as a holomorphic covering
$f_m:\overline{X}_m\to \overline{\CP}^2$, also branched along $C_{m,1}$.
Consider the standard complex conjugation 
$c:\CP^2 \to \CP^2 $, 
in homogeneous coordinates it is given by $c^*(x_i)=\overline{x}_i$. 
The composition 
$\overline{f}_m=c\circ f_m :\overline{X}_m\to \CP^2$ is a 
holomorphic generic covering with branch curve 
$C_{m,2}=c(C_{m,1})$. By construction, we have 
$$\overline{f}_m^{\,*}(\Lambda )= -f_m^{*}(\Lambda )=-m[K]\,  ,$$ 
where $\Lambda \in H^2(\CP^2,\Q)$ is the class of the projective line in 
$\CP^2$.

The set of generic coverings $f$ of $\CP^2$ branched along a cuspidal curve 
$C$ is in one-to-one correspondence with the set of epimorphisms from 
the fundamental group $\pi _1(\CP^2\setminus C)$ to the symmetric groups 
$S_{\deg f}$ (up to inner automorphisms) satisfying some additional 
properties (see \cite{K}). 
By Theorem 3 in \cite{K}, for $C_{m,1}$ (respectively, $C_{m,2}$) 
there exists only one such an epimorphism
$\varphi_m: \pi _1(\CP^2\setminus C_{m,1})\to S_{\deg f_m}$  
(respectively, $\overline{\varphi}_m$). 

Assume that there is an   
isotopy
$F_t: \CP^2 \to \CP^2$ such that 
$F_1 (C_{m,1})=C_{m,2}$.
The epimorphism 
$$
\varphi_{m,t}: \pi _1(\CP^2\times [0,1]\setminus 
\{ (F_t(C_{m,1}),t)\} )\to S_{\deg f_m}
$$ 
defines a  
5-manifold 
$Y$ (with boundary) as a generic covering of 
$f_{m,t}:Y\to
\CP^2\times [0,1]$. This manifold has a natural structure of a 
locally trivial fibration over  
$[0,1]$.  
Note that by Theorem 3 in \cite{K}, 
over $t=0$ the covering coincides 
with $f_m$ and over $t=1$ with $\overline{f}_m$. 
This fibration provides a homeomorphism $h$ of $X$ 
such that $h^*[K]=-[K]$. 
Thus, the isotopy $F_t$ can not exist. (Note that, in fact, $Y$ is smooth
and $h$ can be made smooth everywhere except a finite number of
points.)
 
If we assume that the branch curves
$C_{m,1}$ $C_{m,2}$ have the same braid factorization type, 
then Theorem 7.1 in \cite{KT}  
would imply the existence of 
an 
isotopy
$F_t : \CP^2 \to \CP^2$ such that $F_1 (C_{m,1})=C_{m,2}$. 
This completes the proof of Theorem 2.

To prove the existence of curves satisfying Theorem 1, consider the arrangement of nine
lines $L=L_1\cup \dots \cup L_9$  
in $\CP ^2$ dual to the nine inflection points of a smooth cubic $C$
in the dual plane.
If $C$ is a cubic given by
$x_1^3+x_2^3+x_3^3=0$,
then the lines
$L_1, \dots L_9$ are given by equations
$$
\matrix
L_1=\{x_1-x_3=0 \}, & L_2=\{x_1-\mu^2 x_3=0 \}, &
L_3=\{x_1+\mu x_3=0 \}, \\
L_4=\{x_2-\mu^2x_3=0 \}, & L_5=\{x_2-x_3=0 \}, &
L_6=\{x_2+\mu x_3=0  \},  \\
L_7=\{x_1+\mu x_2=0  \}, & L_8=\{x_1-\mu^2x_2=0 \}, &
L_9=\{ x_1-x_2=0 \}.
\endmatrix
$$
where $\mu =e^{\pi i/3}$.

Consider an affine surface $S\subset \C ^4$ given in $\C^4$ by equations
$$
\align
&w^5_1=l_1l_2l_3l_4^3l_5^3l_9\,,\\
&w_2=l_1l_3l_4^3l_6l_7l_8^2l_9\,,\\
\endalign
$$
where $l_i=0$ are the equations of $L_i$ in the chart $\{ x_3\neq 0\} $. 
Denote by X  
the minimal desingularization of the normalization of the projective closure 
of $S$ in $\CP^4$. 

As is shown in \cite{KK}, $X$ is a rigid surface of general type with 
$K_X^2=333$, the topological Euler characteristic $e(X)=111$, and such that 
$X$ has no anti-holomorphic automorphisms. Its universal covering is the 
complex ball. By 
Mostow rigidity \cite{M},
any homeomorphism $h$ of $X$ is homotopic to a
holomorphic or anti-holomorphic automorphism of $X$.
Thus, $X$ satisfies the conditions of Theorem 2.

Consider a generic covering $f_m: X_m\to \CP ^2$ given by $m$-canonical 
class and let 
$C_{m,1}\subset \CP^2$ be the branch curve of $f_m$. 
Calculations in \cite{K}, page 1155,
give rise to 
\roster
\item"{(1)}" $\deg f_m=333\,m^2$;
\item"{(2)}" $C_{m,1}$ is a plane cuspidal curve of 
$\deg C_{m,1}= 333\, m(3m+1)$;
\item"{(3)}" the geometric genus $g_m$ of $C_{m,1}$ equals
$g_m=333\,(3m+2)(3m+1)/2+1$;
\item"{(4)}" $C_{m,1}$ has $c=111\,(36m^2+27m+5)$ ordinary cusps.
\endroster  
As is proved above,
the curves $C_{m,1}$ and $C_{m,2}=\overline{C}_{m,1}$ give  
a sequence of curves satisfying Theorem 1. 
(Note that other  examples of such curves can be obtained in 
a similar way 
if we take for $X$
a fake projective plane, 
i.e., a surface  
of general type with $p_g=q=0$ and $K^2=9$; 
the universal covering of such surfaces is again
the complex ball and, 
as is shown in \cite{KK}, they have no anti-isomorphisms either.)\qed\enddemo

\end